\documentclass{amsart}

\usepackage{amssymb,latexsym}
\usepackage[latin1]{inputenc} 
\usepackage[dvips]{graphicx,color,psfrag}

\numberwithin{equation}{section}
\numberwithin{figure}{section}
\numberwithin{table}{section}

\theoremstyle{plain}
\newtheorem{prop}{Proposition}[section]
\newtheorem{lemma}[prop]{Lemma}
\newtheorem{theorem}[prop]{Theorem}

\newtheorem{conjecture}[prop]{Conjecture}

\theoremstyle{definition}
\newtheorem{defin}{Definition}[section]

\newenvironment{proofoft}[1]{{\em Proof of Theorem #1. }}{$\Box$ \vspace{1em}}

\newenvironment{proofofprop}[1]{{\em Proof of Proposition #1. }}{$\Box$ \vspace{1em}}
\newenvironment{proofofprop_no_text}{}{$\Box$ \vspace{1em}}

\newcommand{\Zd}{{\mathbb Z}^d}
\newcommand{\T}{{\mathbb T}}
\newcommand{\Td}{{\mathbb T}^d}

\newcommand{\RR}{{\mathbb R}}

\newcommand{\Remark}{\em Remark. \rm}
\newcommand{\Remarks}{\em Remarks. \rm}

\newcommand{\dpt}{{}^{\,\,}\!} 
\newcommand{\dptn}{\,{}^{\!\!}}
\newcommand{\scalefactor}{0.87}
\def\P{{\mathbf P}}

\begin{document}

\title[Stochastic domination for the Ising and fuzzy Potts models]
{Stochastic domination for the Ising and fuzzy Potts models}
\date{\today}

\author[Marcus Warfheimer]{Marcus Warfheimer}
\address{Department of Mathematical Sciences, Chalmers University of Technology and University of Gothenburg, SE--41296 Gothenburg, Sweden}
\email{marcus.warfheimer@gmail.com}
\urladdr{http://www.math.chalmers.se/~warfheim}
\thanks{Research partially supported by the Göran Gustafsson Foundation for Research in Natural Sciences and Medicine.}



\begin{abstract}
We discuss various aspects concerning stochastic domination for the Ising model and the fuzzy Potts model. We begin by considering the Ising model on the homogeneous tree of degree $d$, $\Td$. For given interaction parameters $J_1$, $J_2>0$ and external field $h_1\in\RR$, we compute the smallest external field $\tilde{h}$ such that the plus measure with parameters $J_2$ and $h$ dominates the plus measure with parameters $J_1$ and $h_1$ for all $h\geq\tilde{h}$. Moreover, we discuss continuity of $\tilde{h}$ with respect to the three parameters $J_1$, $J_2$, $h$ and also how the plus measures are stochastically ordered in the interaction parameter for a fixed external field. Next, we consider the fuzzy Potts model and prove that on $\Zd$ the fuzzy Potts measures dominate the same set of product measures while on $\Td$, for certain parameter values, the free and minus fuzzy Potts measures dominate different product measures. For the Ising model, Liggett and Steif proved that on $\Zd$ the plus measures dominate the same set of product measures while on $\T^2$ that statement fails completely except when there is a unique phase.  
\end{abstract}

\maketitle

\section{Introduction and main results}

The concept of stochastic domination has played an important role in probability theory over the last couple of decades, for example in interacting particle systems and statistical mechanics. In \cite{Liggett_Steif}, various results were proved concerning stochastic domination for the Ising model with no external field on $\Zd$ and on the homogeneous binary tree $\T^2$ (i.e.\ the unique infinite tree where each site has $3$ neighbors). As an example, the following distinction between $\Zd$ and $\T^2$ was shown: On $\Zd$, the plus and minus states dominate the same set of product measures, while on $\T^2$ that statement fails completely except in the case when we have a unique phase. In this paper we study stochastic domination for the Ising model in the case of nonzero external field and also for the so called fuzzy Potts model. 

Let $V$ be a finite or countable set and equip the space $\{-1,1\}^V$ with the following natural partial order: For $\eta$, $\eta^\prime\in\{-1,1\}^V$, we write $\eta\leq\eta^\prime$ if $\eta(x)\leq\eta^\prime(x)$ for all $x\in V$. Moreover, whenever we need a topology on $\{-1,1\}^V$ we will use the product topology. We say that a function $f: \{-1,1\}^V\to\RR$ is increasing if $f(\eta)\leq f(\eta^\prime)$ whenever $\eta\leq\eta^\prime$. We will use the following usual definition of stochastic domination.
\begin{defin}[Stochastic domination]
Given a finite or countable set $V$ and probability measures $\mu_1$, $\mu_2$ on $\{-1,1\}^V$, we say that $\mu_2$ dominates $\mu_1$ (written $\mu_1\leq\mu_2$ or $\mu_2\geq\mu_1$) if
\[
\int f \,d\mu_1 \leq \int f \,d\mu_2
\]
for all real-valued, continuous and increasing functions $f$ on $\{-1,1\}^V$. 
\end{defin}
It is well known that a necessary and sufficient condition for two probability measures $\mu_1$, $\mu_2$ to satisfy $\mu_1\leq\mu_2$ is that there exists a coupling measure $\nu$ on $\{-1,1\}^V\times\{-1,1\}^V$ with first and second marginals equal to $\mu_1$ and $\mu_2$ respectively and 
\[
\nu(\,(\eta,\xi):\eta\leq\xi\,)=1.
\]
(For a proof, see for example \cite[p.~72-74]{Liggett85}.) Given any set $S\subseteq \RR$ and a family of probability measures $\{\mu_s\}_{s\in S}$ indexed by $S$, we will say that the map $S\ni s\mapsto \mu_s$ is increasing if $\mu_{s_1}\leq \mu_{s_2}$ whenever $s_1<s_2$.

\subsection{The Ising model}

The ferromagnetic Ising model is a well studied object in both physics and probability theory. For a given infinite, locally finite (i.e.\ each vertex has a finite number of neighbors), connected graph $G=(V,E)$, it is defined from the nearest-neighbor potential
\[
\Phi_A^{J,h}(\eta)=\begin{cases}
-J\eta(x)\eta(y)& \text{if $A=\{x,y\}$, with $\langle x, y\rangle\in E$},\\
-h\eta(x)& \text{if $A=\{x\}$},\\
0& \text{otherwise}
\end{cases}
\]    
where $A\subseteq V$, $\eta\in\{-1,1\}^{V}$, $J>0$, $h\in \RR$ are two parameters called the coupling strength and the external field respectively and $\langle x,y\rangle$ denotes the edge connecting $x$ and $y$. A probability measure $\mu$ on $\{-1,1\}^{V}$ is said to be a Gibbs measure (or sometimes Gibbs state) for the ferromagnetic Ising model with parameters $h\in\RR$ and $J>0$ if it admits conditional probabilities such that for all finite $U\subseteq V$, all $\sigma\in\{-1,1\}^U$ and all $\eta\in\{-1,1\}^{V\setminus U}$
\begin{align*}
&\mu(X(U)=\sigma \,|\,X(V\setminus U)=\eta)\\
&\quad=\frac{1}{Z_{J,h}^{U,\eta}}\exp\Bigg[J\Bigg(\displaystyle\sum_{\langle x,y \rangle\in E, x,y\in U}\sigma(x)\sigma(y)+\sum_{\langle x,y \rangle\in E, x\in U, y\in\partial U}\sigma(x)\eta(y)\Bigg) \\
&\qquad+h\sum_{x\in U}\sigma(x)\Bigg]
\end{align*}
where $Z_{J,h}^{U,\eta}$ is a normalizing constant and 
\[
\partial U=\{\,x \in V\setminus U:\exists y\in U \text{ such that } \langle x,y\rangle\in E\,\}.
\]
For given $J>0$ and $h\in\RR$, we will denote the set of Gibbs measures with parameters $J$ and $h$ by $\mathcal G(J,h)$ and we say that a phase transition occurs if $|\mathcal G(J,h)|>1$, i.e.\ if there exist more than one Gibbs state. (From the general theory described in \cite{Georgii} or \cite{Liggett85}, $\mathcal G(J,h)$ is always nonempty.) At this stage one can ask, for fixed $h\in\RR$, is it the case that the existence of multiple Gibbs states is increasing in $J$? When $h=0$ it is possible from the so called random-cluster representation of the Ising model to show a positive answer to the last question (see \cite{Haggstrom_random_clust_repr} for the case when $G=\Zd$ and \cite{Haggstrom_gen_graphs} for more general $G$). However, when $h\neq 0$ there are graphs where the above monontonicity property no longer holds, see \cite{Schonmann_Tanaka_mon} for an example of a relatively simple such graph. 

Furthermore, still for fixed $J> 0$, $h\in\RR$, standard monotonicity arguments can be used to show that there exists two particular Gibbs states $\mu_h^{J,+}$, $\mu_h^{J,-}$, called the plus and the minus state, which are extreme with respect to the stochastic ordering in the sense that 
\begin{equation}\label{paper3:extreme}
\mu_h^{J,-}\leq \mu\leq\mu_h^{J,+} \quad \text{for any other $\mu\in\mathcal G(J,h)$}. 
\end{equation}
To simplify the notation, we will write $\mu^{J,+}$ for $\mu_0^{J,+}$ and $\mu^{J,-}$ for $\mu_0^{J,-}$. (Of course, most of the things we have defined so far are also highly dependent on the graph $G$, but we suppress that in the notation.) 

In \cite{Liggett_Steif} the authors studied, among other things, stochastic domination between the plus measures $\{\mu^{J,+}\}_{J>0}$ in the case when $G=\T^2$. For example they showed that the map $(0,\infty)\ni J\mapsto \mu^{J,+}$ is increasing when $J>J_c$ and proved the existence of and computed the smallest $J>J_c$ such that $\mu^{J,+}$ dominates $\mu^{J^\prime,+}$ for all $0<J^\prime \leq J_c$. (On $\Zd$, the fact that $\mu^{J_1,+}$ and $\mu^{J_2,+}$ are not stochastically ordered when $J_1\neq J_2$ gives that such a $J$ does not even exist in that case.) Our first result deals with the following question: Given $J_1$, $J_2>0$, $h_1\in\RR$, can we find the smallest external field $\tilde{h}=\tilde{h}(J_1,J_2,h_1)$ with the property that $\mu_h^{J_2,+}$ dominates $\mu_{h_1}^{J_1,+}$ for all $h\geq\tilde{h}$? To clarify the question a bit more, note that an easy application of Holley's theorem (see \cite{Georgii_Haggstrom_Maes_geom}) tells us that for fixed $J>0$, the map $\RR\ni h\mapsto\mu_h^{J,+} $ is increasing. Hence, for given $J_1$, $J_2$ and $h_1$ as above the set
\[
\{\,h\in\RR:\mu_h^{J_2,+}\geq \mu_{h_1}^{J_1,\pm}\,\}
\]
is an infinite interval and we want to find the left endpoint of that interval (possibly $-\infty$ or $+\infty$ at this stage). For a general graph not much can be said, but we have the following easy bounds on $\tilde{h}$ when $G$ is of bounded degree.
\begin{prop}\label{paper3:gen_graph}
Consider the Ising model on a general graph $G=(V,E)$ of bounded degree. Define
\[
\tilde{h}=\tilde{h}(J_1,J_2,h_1)=\inf\{\,h\in\RR:\,\mu_h^{J_2,+}\geq\mu_{h_1}^{J_1,+}\,\}.
\]
Then 
\[
h_1-N(J_1+J_2)\leq \tilde{h}\leq h_1+N|J_1-J_2|,
\]
where $N=\displaystyle\sup_{x\in V}N_x$ and $N_x$ is the number of neighbors of the site $x\in V$. 
\end{prop}

For the Ising model, we will now consider the case when $G=\T^d$, the homogeneous $d$-ary tree, defined as the unique infinite tree where each site has exactly $d+1\geq 3$ neighbors. The parameter $d$ is fixed in all that we will do and so we suppress that in the notation. For this particular graph it is well known that for given $h\in\RR$, the existence of multiple Gibbs states is increasing in $J$ and so as a consequence there exists a critical value $J_c(h)\in [0,\infty]$ such that when $J<J_c(h)$ we have a unique Gibbs state whereas for $J>J_c(h)$ there are more than one Gibbs states. In fact, much more can be shown in this case. As an example it is possible to derive an explicit expression for the phase transition region 
\[
\{\,(J,h)\in \RR^2:\,|\mathcal G(J,h)|>1\,\},
\]
in particular one can see that $J_c(h)\in (0,\infty)$ for all $h\in\RR$. Moreover,
\[
J_c:=J_c(0)=\text{arccoth }\dptn d=\frac{1}{2}\log\frac{d+1}{d-1},
\]
see \cite{Georgii} for more details. (Here and in the sequel, $:=$ will mean definition.)

To state our results for the Ising model on $\T^d$, we need to recall some more facts, all of which can be found in \cite[p.~247-255]{Georgii}. To begin, we just state what we need very briefly and later on we will give some more details. Given $J>0$ and $h\in\RR$, there is a one-to-one correspondence $t\mapsto \mu$ between the real solutions of a certain equation (see \eqref{paper3:fixpunkt} and the function $\phi_J$ in \eqref{paper3:phi_J} below) and the completely homogeneous Markov chains in $\mathcal G(J,h)$ (to be defined in Section~\ref{paper3:proofs}). Let $t_\pm(J,h)$ denote the real numbers which correspond to the plus and minus measure respectively. (It is easy to see that the plus and minus states are completely homogeneous Markov chains, see Section~\ref{paper3:proofs}.) We will write $t_\pm(J)$ instead of $t_\pm(J,0)$. Furthermore, let 
\[
h^*(J)=\displaystyle\max_{t\geq 0}\big(d\dpt\phi_J(t)-t\big) 
\]
and denote by $t^*(J)$ the $t\geq 0$ where the function $t\mapsto d\dpt\phi_J(t)-t$ attains its unique maximum. In \cite{Georgii}, explicit expressions for both $h^*$ and $t^*$ are derived:
\begin{align*}
&h^*(J)=\begin{cases}
0 \quad &\text{if $J\leq J_c$} \\
d\dpt\text{arctanh}\left(\frac{d\tanh(J)-1}{d\coth(J)-1}\right)^{1/2}- \text{arctanh}\left(\frac{d-\coth(J)}{d-\tanh(J)}\right)^{1/2} \quad &\text{if $J>J_c$} 
\end{cases} \\ 
&t^*(J)=\begin{cases}
0 \quad &\text{if $J\leq J_c$} \\
\text{arctanh}\left(\frac{d-\coth(J)}{d-\tanh(J)}\right)^{1/2} \quad &\text{if $J>J_c$}\end{cases} 
\end{align*}

In particular one can see that both $h^*$ and $t^*$ are continuous functions of $J$ and by computing derivatives one can show that they are strictly increasing for $J>J_c$.
\begin{figure}[!h]
\begin{center}
\scalebox{\scalefactor}{
%
%
\begin{psfrags}%
\psfragscanon%
%
\psfrag{s05}[t][t]{\color[rgb]{0,0,0}\setlength{\tabcolsep}{0pt}\begin{tabular}{c}Interaction parameter $J$\end{tabular}}%
\psfrag{s06}[b][b]{\color[rgb]{0,0,0}\setlength{\tabcolsep}{0pt}\begin{tabular}{c}$h^*(J)$\end{tabular}}%
\psfrag{s07}[t][t]{\color[rgb]{0,0,0}\setlength{\tabcolsep}{0pt}\begin{tabular}{c}Interaction parameter $J$\end{tabular}}%
\psfrag{s08}[b][b]{\color[rgb]{0,0,0}\setlength{\tabcolsep}{0pt}\begin{tabular}{c}$t^*(J)$\end{tabular}}%
%
\psfrag{x01}[t][t]{0}%
\psfrag{x02}[t][t]{0.2}%
\psfrag{x03}[t][t]{0.4}%
\psfrag{x04}[t][t]{0.6}%
\psfrag{x05}[t][t]{0.8}%
\psfrag{x06}[t][t]{1}%
\psfrag{x07}[t][t]{0}%
\psfrag{x08}[t][t]{0.2}%
\psfrag{x09}[t][t]{0.4}%
\psfrag{x10}[t][t]{0.6}%
\psfrag{x11}[t][t]{0.8}%
\psfrag{x12}[t][t]{1}%
%
\psfrag{v01}[r][r]{-1}%
\psfrag{v02}[r][r]{0}%
\psfrag{v03}[r][r]{1}%
\psfrag{v04}[r][r]{2}%
\psfrag{v05}[r][r]{-1}%
\psfrag{v06}[r][r]{0}%
\psfrag{v07}[r][r]{1}%
\psfrag{v08}[r][r]{2}%
%
\resizebox{12cm}{!}{\includegraphics{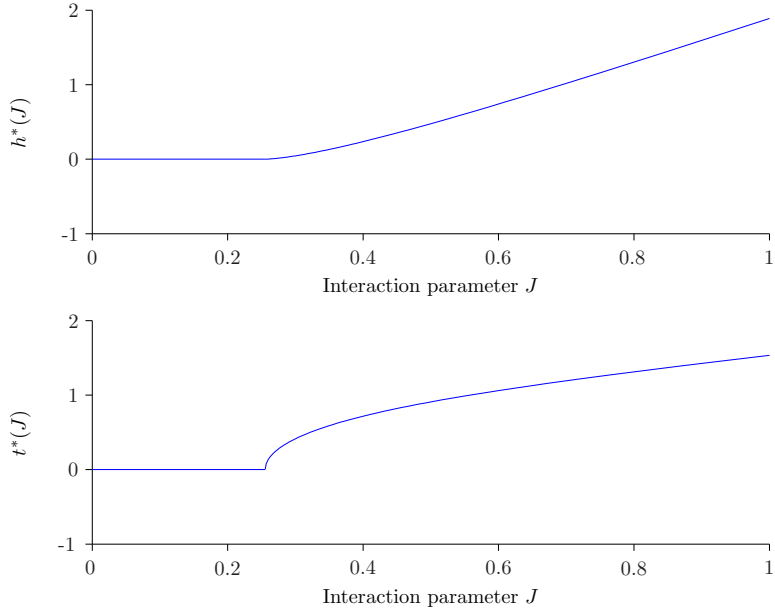}}%
\end{psfrags}%
%
}
\caption{The functions $h^*$ and $t^*$ in the case when $d=4$.}
\label{paper3:hstar}
\end{center}
\end{figure}
\begin{theorem}\label{paper3:prop_dom_tree}
Consider the Ising model on $\T^d$ and let $J_1$, $J_2>0$, $h_1\in\RR$ be given. Define
\begin{align*}
f_{\pm}(J_1,J_2,h_1)&=\inf\{\,h\in\RR:\mu_h^{J_2,+}\geq \mu_{h_1}^{J_1,\pm}\,\} \\
g_{\pm}(J_1,J_2,h_1)&=\inf\{\,h\in\RR:\mu_h^{J_2,-}\geq \mu_{h_1}^{J_1,\pm}\,\}
\end{align*}
and denote $\tau_\pm=\tau_{\pm}(J_1,J_2,h_1)=t_{\pm}(J_1,h_1)+|J_1-J_2|$. Then the following holds:
\begin{equation}\label{paper3:prop_dom_tree_ekv1}
f_{\pm}(J_1,J_2,h_1)=\begin{cases}
-h^*(J_2) &\text{if } t_-(J_2,-h^*(J_2))\leq \tau_\pm<t^*(J_2) \\
\tau_\pm-d\dpt\phi_{J_2}(\tau_\pm) &\text{if } \tau_\pm\geq t^*(J_2) \text{ or } \tau_\pm<t_-(J_2,-h^*(J_2))
\end{cases}
\end{equation}
\begin{equation}\label{paper3:prop_dom_tree_ekv2}
g_{\pm}(J_1,J_2,h_1)=\begin{cases}
h^*(J_2) &\text{if } -t^*(J_2)< \tau_\pm\leq t_+(J_2,h^*(J_2)) \\
\tau_\pm-d\dpt\phi_{J_2}(\tau_\pm) &\text{if } \tau_\pm\leq -t^*(J_2) \text{ or } \tau_\pm>t_+(J_2,h^*(J_2))
\end{cases}
\end{equation}
\end{theorem}

\begin{figure}[!h]
\begin{center}
\scalebox{\scalefactor}{
%
%
\begin{psfrags}%
\psfragscanon%
%
\psfrag{s07}[t][t]{\color[rgb]{0,0,0}\setlength{\tabcolsep}{0pt}\begin{tabular}{c}$t$\end{tabular}}%
\psfrag{s08}[l][l]{\color[rgb]{0,0,0}\setlength{\tabcolsep}{0pt}\begin{tabular}{l}$t\mapsto\psi(J_2,t)$\end{tabular}}%
\psfrag{s09}[t][t]{\color[rgb]{0,0,0}\setlength{\tabcolsep}{0pt}\begin{tabular}{c}$t$\end{tabular}}%
\psfrag{s10}[l][l]{\color[rgb]{0,0,0}\setlength{\tabcolsep}{0pt}\begin{tabular}{l}$t\mapsto\theta(J_2,t)$\end{tabular}}%
%
\psfrag{x01}[t][t]{-20}%
\psfrag{x02}[t][t]{-10}%
\psfrag{x03}[t][t]{0}%
\psfrag{x04}[t][t]{10}%
\psfrag{x05}[t][t]{20}%
\psfrag{x06}[t][t]{-20}%
\psfrag{x07}[t][t]{-10}%
\psfrag{x08}[t][t]{0}%
\psfrag{x09}[t][t]{10}%
\psfrag{x10}[t][t]{20}%
%
\psfrag{v01}[r][r]{-10}%
\psfrag{v02}[r][r]{0}%
\psfrag{v03}[r][r]{10}%
\psfrag{v04}[r][r]{-10}%
\psfrag{v05}[r][r]{0}%
\psfrag{v06}[r][r]{10}%
%
\resizebox{12cm}{!}{\includegraphics{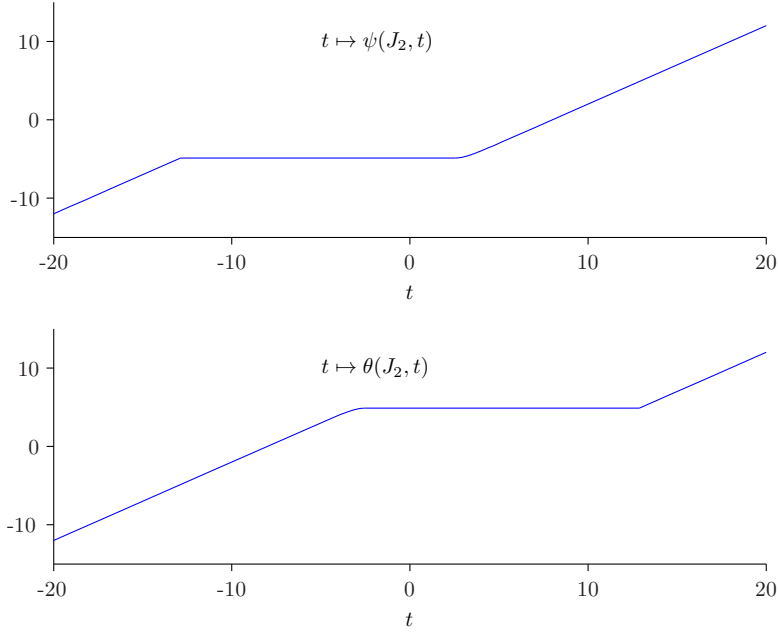}}%
\end{psfrags}%
%
}
\caption{The functions $t\mapsto\psi(J_2,t)$ and $t\mapsto \theta(J_2,t)$ in the case when $J_2=2$ and $d=4$.}
\label{paper3:fplusfig}
\end{center}
\end{figure}

\noindent\Remarks 
\begin{itemize}
\item[(i)] Note that if $0<J_2\leq J_c$, then $h^*(J_2)=0$ and 
\[
t_-(J_2,-h^*(J_2))=t^*(J_2)=t_+(J_2,h^*(J_2))=0 
\]
and hence the first interval disappears in the formulas and we simply get
\[
\begin{split}
f_{\pm}(J_1,J_2,h_1)&=g_{\pm}(J_1,J_2,h_1) \\
&=\tau_{\pm}(J_1,J_2,h_1)-d\phi_{J_2}(\tau_{\pm}(J_1,J_2,h_1)).
\end{split}
\]
\item[(ii)] By looking at the formulas \eqref{paper3:prop_dom_tree_ekv1} and \eqref{paper3:prop_dom_tree_ekv2}, we see that there are functions $\psi$, $\theta:(0,\infty)\times \RR\mapsto \RR$ such that 
\[
\begin{split}
f_{\pm}(J_1,J_2,h_1)&=\psi(J_2,\tau_{\pm}(J_1,J_2,h_1)) \quad \text{and} \\
g_{\pm}(J_1,J_2,h_1)&=\theta(J_2,\tau_{\pm}(J_1,J_2,h_1)). 
\end{split}
\]
(Of course, $\psi(J_2,t)$ and $\theta(J_2,t)$ are just \eqref{paper3:prop_dom_tree_ekv1} and \eqref{paper3:prop_dom_tree_ekv2} with $t$ instead of $\tau_\pm$.) It is easy to check that for fixed $J_2>0$, the maps $t\mapsto \psi(J_2,t)$ and $t\mapsto \theta(J_2,t)$ are continuous. A picture of these functions when $J_2=2$, $d=4$ can be seen in Figure \ref{paper3:fplusfig}. 
\item[(iii)] It is not hard to see by direct computations that $f_+$ satisfies the bounds in Proposition \ref{paper3:gen_graph}. We will indicate how this can be done after the proof of Theorem \ref{paper3:prop_dom_tree}.
\item[(iv)] We will see in the proof that if 
\[
t_-(J_2,-h^*(J_2))\leq \tau_\pm(J_1,J_2,h_1)<t^*(J_2), 
\]
then 
\[
\{\,h\in\RR:\mu_h^{J_2,+}\geq \mu_{h_1}^{J_1,\pm}\,\}=[-h^*(J_2),\infty),
\] 
and if $-t^*(J_2)< \tau_\pm(J_1,J_2,h_1)\leq t_+(J_2,h^*(J_2))$, then 
\[
\{\,h\in\RR:\mu_h^{J_2,-}\geq \mu_{h_1}^{J_1,\pm}\,\}=(h^*(J_2),\infty).
\]
Hence in the first case the left endpoint belongs to the interval, while in the second case it does not.
\end{itemize}

\noindent Our next proposition deals with continuity properties of $f_{\pm}$ and $g_{\pm}$ with respect to the parameters $J_1$, $J_2$ and $h_1$. We will only discuss the function $f_+$, the other ones can be treated in a similar fashion.
\begin{prop}\label{paper3:prop_dom_tree2}
Consider the Ising model on $\T^d$ and recall the notation from Theorem \ref{paper3:prop_dom_tree}. Let
\begin{align*}
a=a(J_1,J_2)&=t_-(J_1,-h^*(J_1))+|J_1-J_2| \\
b=b(J_1,J_2)&=t_+(J_1,-h^*(J_1))+|J_1-J_2|
\end{align*}
\begin{itemize}
\item[a)] Given $J_1$, $J_2>0$, the map $\RR\ni h_1\mapsto f_+(J_1,J_2,h_1)$ is continuous except possibly at $-h^*(J_1)$ depending on $J_1$ and $J_2$ in the following way:
\begin{itemize}
\item[] If $J_1\leq J_c$ or $J_1=J_2$ then it is continuous at $-h^*(J_1)$.
\item[] If $J_1>J_c$ and $0<J_2\leq J_c$ then it is discontinuous at $-h^*(J_1)$.
\item[] If $J_1$, $J_2>J_c$, $J_1\neq J_2$ then it is discontinuous except when
\[
t_-(J_2,-h^*(J_2))\leq a<t^*(J_2) \;\, \text{and} \;\, t_-(J_2,-h^*(J_2))\leq b \leq t^*(J_2).
\]
\end{itemize}
\item[b)] Given $J_2>0$, $h_1\in\RR$, the map $(0,\infty)\ni J_1\mapsto f_+(J_1,J_2,h_1)$ is continuous at $J_1$ if $0<J_1\leq J_c$ or $J_1>J_c$ and $h_1\neq -h^*(J_1)$. In the case when $h_1=-h^*(J_1)$ it is discontinuous at $J_1$ except when 
\[
t_-(J_2,-h^*(J_2))\leq a<t^*(J_2) \;\, \text{and} \;\, t_-(J_2,-h^*(J_2))\leq b \leq t^*(J_2).
\]
\item[c)] Given $J_1>0$, $h_1\in\RR$, the map $(0,\infty)\ni J_2\mapsto f_+(J_1,J_2,h_1)$ is continuous for all $J_2>0$. 
\end{itemize}
\end{prop}

We conclude this section with a result about how the measures $\{\mu_h^{J,+}\}_{J>0}$ are ordered with respect to $J$ for fixed $h\in\RR$.
\begin{prop}\label{paper3:prop_dom_tree3}
Consider the Ising model on $\T^d$. The map $(0,\infty)\ni J\mapsto \mu_h^{J,+}$ is increasing in the following cases: $a)$ $h\geq 0$ and $J\geq J_c$, $b)$ $h<0$ and $h^*(J)>-h$.
\end{prop}

\subsection{The fuzzy Potts model}

Next, we consider the so called fuzzy Potts model. To define the model, we first need to define the perhaps more familiar Potts model. Let $G=(V,E)$ be an infinite locally finite graph and suppose that $q\geq 3$ is an integer. Let $U$ be a finite subset of $V$ and consider the finite graph $H$ with vertex set $U$ and edge set consisting of those edges $\langle x,y\rangle\in E$ with $x,y \in U$. In this way, we say that the graph $H$ is induced by $U$. The finite volume Gibbs measure for the $q$-state Potts model at inverse temperature $J\geq 0$ with free boundary condition is defined to be the probability measure $\pi_{q,J}^H$ on $\{1,2,\dots,q\}^U$ which to each element $\sigma$ assigns probability
\[
\pi_{q,J}^H(\sigma)=\frac{1}{Z_{q,J}^H}\exp\bigg(2J\sum_{\langle x,y \rangle\in E, x,y\in U}I_{\{\sigma(x)=\sigma(y)\}}\bigg),
\] 
where $Z_{q,J}^H$ is a normalizing constant. 

Now, suppose $r\in\{1,\dots,q-1\}$ and pick a $\pi_{q,J}^H\,$-$\,$distributed object $X$ and for $x\in U$ let
\begin{equation}\label{paper3:fuzzy_def}
Y(x)=\begin{cases}
-1& \text{ if } X(x)\in\{1,\dots,r\} \\
\:\:\dpt\dpt\dpt 1 & \text{ if } X(x)\in\{r+1,\dots,q\}.
\end{cases}
\end{equation}
We write $\nu_{q,J,r}^H$ for the resulting probability measure on $\{-1,1\}^U$ and call it the finite volume fuzzy Potts measure on $H$ with free boundary condition and parameters $q$, $J$ and $r$.

We also need to consider the case when we have a boundary condition. For finite $U\subseteq V$, consider the graph $H$ induced by the vertex set $U\cup\partial U$ and let $\eta\in\{1,\dots,q\}^{V\setminus U}$. The finite volume Gibbs measure for the $q$-state Potts model at inverse temperature $J\geq 0$ with boundary condition $\eta$ is defined to be the probability measure on $\{1,\dots,q\}^{U}$ which to each element assigns probability
\[
\begin{split}
\pi_{q,J}^{H,\eta}(\sigma)&=\frac{1}{Z_{q,J}^{H,\eta}}\exp\Bigg(2J\displaystyle\sum_{\langle x,y \rangle\in E, x,y\in U}I_{\{\sigma(x)=\sigma(y)\}} \\
&\quad+2J\sum_{\langle x,y \rangle\in E, x\in U, y\in\partial U}I_{\{\sigma(x)=\eta(y)\}}\Bigg),
\end{split}
\]
where $Z_{q,J}^{H,\eta}$ is a normalizing constant. In the case when $\eta\equiv i$ for some $i\in \{1,\dots,q\}$, we replace $\eta$ with $i$ in the notation.

Furthermore, we introduce the notion of infinite volume Gibbs measure for the Potts model. A probability measure $\mu$ on $\{1,\dots,q\}^V$ is said to be an infinite volume Gibbs measure for the $q$-state Potts model on $G$ at inverse temperature $J\geq 0$, if it admits conditional probabilities such that for all finite $U\subseteq V$, all $\sigma\in\{1,\dots,q\}^U$ and all $\eta\in\{1,\dots,q\}^{V\setminus U}$
\[
\mu(X(U)=\sigma \,|\,X(V\setminus U)=\eta)=\pi_{q,J}^{H,\eta}(\sigma)
\]
where $H$ is the graph induced by $U\cup\partial U$. Let $\{V_n\}_{n\geq 1}$ be a sequence of finite subsets of $V$ such that $V_n\subseteq V_{n+1}$ for all $n$, $V=\bigcup_{n\geq 1}V_n$ and for each $n$, denote by $G_n$ the induced graph by $V_n\cup\partial V_n$. Furthermore, for each $i\in \{1,\dots,q\}$, extend $\pi_{q,J}^{G_n,i}$ (and use the same notation for the extension) to a probability measure on $\{1,\dots,q\}^V$ by assigning with probability one the spin value $i$ outside $V_n$. It is well known (and independent of the sequence $\{V_n\}$) that there for each spin $i\in \{1,\dots,q\}$ exists a infinite volume Gibbs measure $\pi^{G,i}_{q,J}$ which is the weak limit as $n\to\infty$ of the corresponding measures $\pi_{q,J}^{G_n,i}$. Moreover, there exists another infinite volume Gibbs measure denoted $\pi_{q,J}^{G,0}$ which is the limit of $\pi_{q,J}^{G_n}$ in the sense that the probabilities on cylinder sets converge. The existence of the above limits as well as the independence of the choice of the sequence $\{V_n\}$ when constructing them follows from the work of Aizenman et al.\ \cite{Aizenman_et_al_conv_Potts}.

Given the infinite volume Gibbs measures $\{\pi^{G,i}_{q,J}\}_{i\in\{0,\dots,q\}}$, we define the corresponding infinite volume fuzzy Potts measures $\{\nu^{G,i}_{q,J,r}\}_{i\in\{0,\dots,q\}}$ using \eqref{paper3:fuzzy_def}.

In words, the fuzzy Potts model can be thought of arising from the ordinary $q$-state Potts model by looking at a pair of glasses that prevents from distinguishing some of the spin values. From this point of view, the fuzzy Potts model is one of the most basic examples of a so called hidden Markov field \cite{Kunsch_hidden_Markov}. For earlier work on the fuzzy Potts model, see for example \cite{Haggstrom_Gibsian_Zd,Haggstrom_Gibsian_tree,Kahn_frac_fuzzy,Maes_Velde,Haggstrom_pos_corr}.

Given a finite or countable set $V$ and $p\in [0,1]$, let $\gamma_p$ denote the product measure on $\{-1,1\}^V$ with $\gamma_p(\eta:\,\eta(x)=1)=p$ for all $x\in V$. In \cite{Liggett_Steif} the authors proved the following results for the Ising model. (The second result was originally proved for $d=2$ only but it trivially extends to all $d\geq 2$.) 
\begin{prop}[Liggett, Steif]\label{paper3:LS1}
Fix an integer $d\geq 2$ and consider the Ising model on $\Zd$ with parameters $J>0$ and $h=0$. Then for any $p\in [0,1]$, $\mu^{J,+}\geq \gamma_p$ if and only if $\mu^{J,-}\geq \gamma_p$.
\end{prop}
\begin{prop}[Liggett, Steif]\label{paper3:LS2}
Let $d\geq 2$ be a given integer and consider the Ising model on $\Td$ with paramteters $J>0$ and $h=0$. Moreover, let $\mu^{J,f}$ denote the Gibbs state obtained by using free boundary conditions. If $\mu^{J,+}\neq \mu^{J,-}$, then there exist $0<p^\prime<p$ such that $\mu^{J,+}$ dominates $\gamma_p$ but $\mu^{J,f}$ does not dominate $\gamma_p$ and $\mu^{J,f}$ dominates $\gamma_{p^\prime}$ but $\mu^{J,-}$ does not dominate $\gamma_{p^\prime}$.
\end{prop}
In words, on $\Zd$ the plus and minus state dominate the same set of product measures while on $\Td$ that is not the case except when the we have a unique phase. 

To state our next results we will take a closer look at the construction of the infinite volume fuzzy Potts measures when $G=\Zd$ or $G=\Td$. In those cases it follows from symmetry that $\nu_{q,J,r}^{G,i}=\nu_{q,J,r}^{G,j}$ if $i,j\in\{1,\dots,r\}$ or $i,j\in\{r+1,\dots,q\}$, i.e.\ when the Potts spins $i,j$ map to the same fuzzy spin. For that reason, we let $\nu_{q,J,r}^{G,-}:=\nu_{q,J,r}^{G,1}$ and $\nu_{q,J,r}^{G,+}:=\nu_{q,J,r}^{G,q}$ when $G=\Zd$ or $\Td$. (Of course, we stick to our earlier notation of $\nu_{q,J,r}^{G,0}$.) Our first result is a generalization of Proposition \ref{paper3:LS1} to the fuzzy Potts model. 
\begin{prop}\label{paper3:fuzzy_prod_Zd}
Let $d\geq 2$ be a given integer and consider the fuzzy Potts model on $\Zd$ with parameters $q\geq 3$, $J>0$ and $r\in\{1,\dots,q-1\}$. Then for any $k,l\in\{0,-,+\}$ and $p\in [0,1]$, $\nu^{\Zd,k}_{q,J,r}\geq\gamma_p$ if and only if $\nu^{\Zd,l}_{q,J,r}\geq\gamma_p$. 
\end{prop}
In the same way as for the Ising model, we believe that Proposition \ref{paper3:fuzzy_prod_Zd} fails completely on $\Td$ except when we have a unique phase in the Potts model. Our last result is in that direction. 
\begin{prop}\label{paper3:fuzzy_prod_Td}
Let $d\geq 2$ be a given integer and consider the fuzzy Potts model on $\Td$ with parameters $q\geq 3$, $J>0$ and $r\in\{1,\dots,q-1\}$ where $e^{2J}\geq q-2$. If the underlying Gibbs measures for the Potts model satisfy $\pi^{\Td,1}_{q,J}\neq \pi^{\Td,0}_{q,J}$, then there exists $0<p<1$ such that $\nu^{\Td,0}_{q,J,r}$ dominates $\gamma_p$ but $\nu^{\Td,-}_{q,J,r}$ does not dominate $\gamma_p$. 
\end{prop}

\section{Proofs}\label{paper3:proofs}

We start to recall some facts from \cite{Georgii} concerning the notion of completely homogeneous Markov chains on $\T^d$. Denote the vertex set and the edge set of $\T^d$ with $V(\T^d)$ and $E(\T^d)$ respectively. Given a directed edge $\langle x,y\rangle\in E(\T^d)$ define the ``past'' sites by 
\[
]-\infty,\langle x,y\rangle[=\{\,z\in V(\Td): \text{$z$ is closer to $x$ than to $y$}\,\}.
\] 
For $A\subseteq V(\T^d)$ denote by $\mathcal F_A$ the $\sigma$-algebra generated by the spins in $A$. A probability measure $\mu$ on $\{-1,1\}^{V(\T^d)}$ is called a Markov chain if 
\[
\mu(\,\eta(y)=1\,|\,\mathcal F_{]-\infty,\langle x,y\rangle[}\,)=\mu(\,\eta(y)=1\,|\,\mathcal F_{\{x\}}\,) \quad \mu\text{-a.s.}
\]
for all $\langle x,y\rangle\in E(\T^d)$. Furthermore, a Markov chain $\mu$ is called completely homogeneous with transition matrix $P=\{\,P(i,j):\,i,j\in\{\ -1,1\}\,\}$ if 
\begin{equation}\label{paper3:compl_hom_prop}
\mu(\,\eta(y)=u\,|\,\mathcal F_{\{x\}}\,)=P(\eta(x),u) \quad \mu\text{-a.s.}
\end{equation}
for all $\langle x,y\rangle\in E(\T^d)$ and $u\in\{-1,1\}$. Observe that such a $P$ necessarily is a stochastic matrix and if it in addition is irreducible denote its stationary distribution by $\nu$. In that situation, we get for each finite connected set $C\subseteq V(\T^d)$, $z\in C$ and $\xi\in\{-1,1\}^{C}$ that
 \[
\mu(\eta=\xi)=\nu(\xi(z))\prod_{\langle x,y\rangle\in D}P(\xi(x),\xi(y))
\] 
where $D$ is the set of directed edges $\langle x,y\rangle$, where $x,y\in C$ and $x$ is closer to $z$ than $y$ is. In particular, it follows that every completely homogeneous Markov chain which arise from an irreducible stochastic matrix is invariant under all graph automorphisms.

Next, we give a short summary from \cite{Georgii} of the Ising model on $\Td$. For $J>0$, define
\begin{equation}\label{paper3:phi_J}
\phi_J(t)=\frac{1}{2}\log\frac{\cosh(t+J)}{\cosh(t-J)},\quad t\in\RR.
\end{equation}
The function $\phi_J$ is trivially seen to be odd. Moreover, $\phi_J$ is concave on $[0,\infty)$, increasing and bounded. (In fact, $\phi_J(t)\to J$ as $t\to\infty$.) Furthermore, there is a one-to-one correspondence $t\mapsto\mu_t$ between the completely homogeneous Markov chains in $\mathcal G(J,h)$ and the numbers $t\in\RR$ satisfying the equation
\begin{equation}\label{paper3:fixpunkt}
t=h+d\,\phi_J(t).
\end{equation} 
In addition, the transition matrix $P_t$ of $\mu_t$ is given by
\begin{equation}\label{paper3:Pt}
\begin{pmatrix}
P_t(-1,-1)& P_t(-1,1)\\
P_t(1,-1)& P_t(1,1)
\end{pmatrix}=
\begin{pmatrix}
\frac{e^{J-t}}{2\cosh(J-t)}& \frac{e^{t-J}}{2\cosh(J-t)}\\
\frac{e^{-J-t}}{2\cosh(J+t)}& \frac{e^{J+t}}{2\cosh(J+t)}\\
\end{pmatrix}.
\end{equation}
\begin{figure}[!h]
\begin{center}
\scalebox{\scalefactor}{
%
%
\begin{psfrags}%
\psfragscanon%
%
\psfrag{s04}[t][t]{\color[rgb]{0,0,0}\setlength{\tabcolsep}{0pt}\begin{tabular}{c}$t$\end{tabular}}%
\psfrag{s05}[l][l]{\color[rgb]{0,0,0}\setlength{\tabcolsep}{0pt}\begin{tabular}{l}$t\mapsto h+d\phi_J(t)$\end{tabular}}%
%
\psfrag{x01}[t][t]{-20}%
\psfrag{x02}[t][t]{-10}%
\psfrag{x03}[t][t]{0}%
\psfrag{x04}[t][t]{10}%
\psfrag{x05}[t][t]{20}%
%
\psfrag{v01}[r][r]{-20}%
\psfrag{v02}[r][r]{-10}%
\psfrag{v03}[r][r]{0}%
\psfrag{v04}[r][r]{10}%
\psfrag{v05}[r][r]{20}%
%
\resizebox{12cm}{!}{\includegraphics{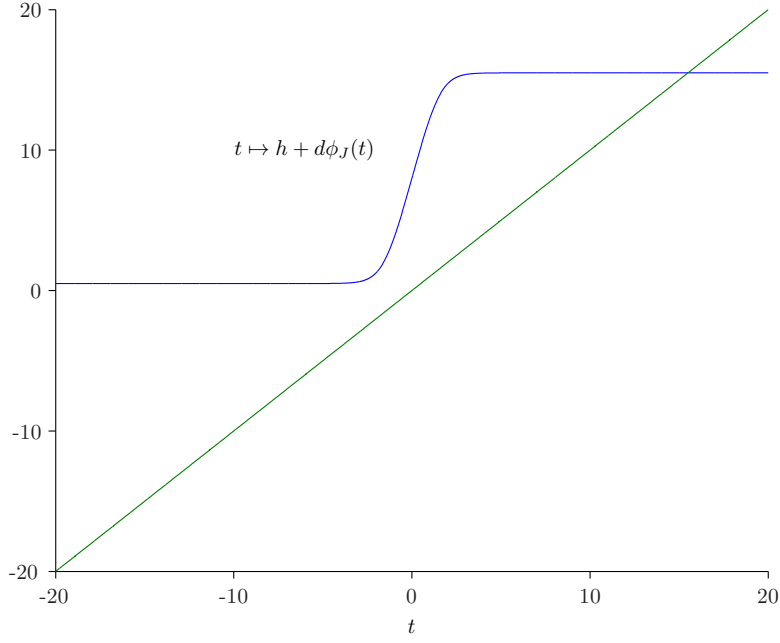}}%
\end{psfrags}%
%
}
\caption{A picture of the fixed point equation \eqref{paper3:fixpunkt} when $d=5$, $h=8$ and $J=3/2$. In this particular case we have a unique solution.}
\label{paper3:fixedplot1}
\end{center}
\end{figure}
\begin{figure}[!h]
\begin{center}
\scalebox{\scalefactor}{
%
%
\begin{psfrags}%
\psfragscanon%
%
\psfrag{s04}[t][t]{\color[rgb]{0,0,0}\setlength{\tabcolsep}{0pt}\begin{tabular}{c}$t$\end{tabular}}%
\psfrag{s05}[l][l]{\color[rgb]{0,0,0}\setlength{\tabcolsep}{0pt}\begin{tabular}{l}$\big(t^*(J),\phi_J(t^*(J))\big)$\end{tabular}}%
%
\psfrag{x01}[t][t]{-10}%
\psfrag{x02}[t][t]{-5}%
\psfrag{x03}[t][t]{0}%
\psfrag{x04}[t][t]{5}%
\psfrag{x05}[t][t]{10}%
%
\psfrag{v01}[r][r]{-10}%
\psfrag{v02}[r][r]{-5}%
\psfrag{v03}[r][r]{0}%
\psfrag{v04}[r][r]{5}%
\psfrag{v05}[r][r]{10}%
%
\resizebox{12cm}{!}{\includegraphics{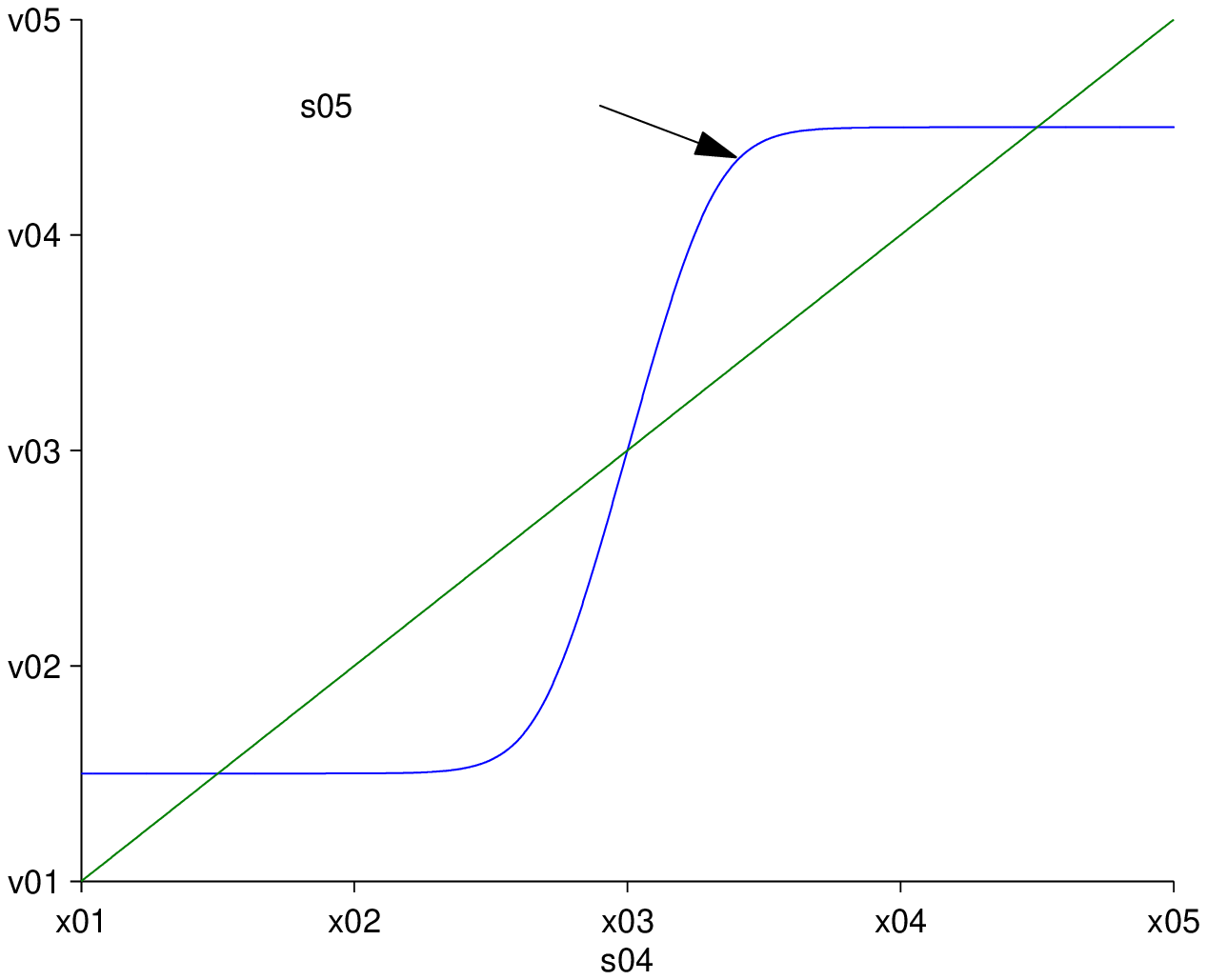}}%
\end{psfrags}%
%
}
\caption{A picture of the fixed point equation \eqref{paper3:fixpunkt} when $d=5$, $h=0$ and $J=3/2$.}
\label{paper3:fixedplot2}
\end{center}
\end{figure}
Given $h\in\RR$ and $J>0$ the fixed point equation \eqref{paper3:fixpunkt} has one, two or three solutions. In fact Lemma \ref{paper3:fixedpoint_Georgii} below tells us exactly when the different situations occur. The largest solution, denoted $t_+(J,h)$, corresponds to the plus measure $\mu_h^{J,+}$ and the smallest, denoted $t_-(J,h)$, to the minus measure $\mu_h^{J,-}$. To see why the last statement is true, let $\mu_\pm=\mu_{t_\pm(J,h)}$ and note that Lemma \ref{paper3:LS3} from Section~\ref{paper3:proof_dom_tree} implies that $\mu_-\leq\mu\leq\mu_+$ for any $\mu\in\mathcal G(J,h)$ which is also a completely homogeneous Markov chain on $\Td$. Moreover, equation \eqref{paper3:extreme} implies that $\mu_h^{J,-}\leq\mu_\pm\leq\mu_h^{J,+}$ and so $\mu_\pm=\mu_h^{J,\pm}$ will follow if $\mu_h^{J,\pm}$ are completely homogeneous Markov chains. To see that, note that equation \eqref{paper3:extreme} also implies that $\mu_h^{J,\pm}$ are extremal in $\mathcal G(J,h)$ which in turn (see Theorem 12.6 in \cite{Georgii}) gives us that they are Markov chains on $\Td$. Finally, from the fact that $\mu_h^{J,\pm}$ are invariant under all graph automorphisms on $\Td$, we obtain the completely homogeneous property \eqref{paper3:compl_hom_prop}.

\begin{lemma}[Georgii]\label{paper3:fixedpoint_Georgii}
The fixed point equation \eqref{paper3:fixpunkt} has
\begin{itemize}
\item[a)] a unique solution when $|h|>h^*(J)$ or $h=h^*(J)=0$,
\item[b)] two distinct solutions $t_-(J,h)<t_+(J,h)$ when $|h|=h^*(J)>0$,
\item[c)] three distinct solutions $t_-(J,h)<t_0(J,h)<t_+(J,h)$ when $|h|<h^*(J)$.
\end{itemize}
\end{lemma} 

\subsection{Proof of Proposition \ref{paper3:gen_graph}}

\begin{proofofprop_no_text}
For the upper bound, just invoke Proposition $4.16$ in \cite{Georgii_Haggstrom_Maes_geom} which gives us that $\mu_h^{J_2,+}\geq \mu_{h_1}^{J_1,+}$ if $h\geq h_1+N|J_1-J_2|$.

For the lower bound, we argue by contradiction as follows. Assume 
\[
\tilde{h}<h_1-N(J_1+J_2)
\] 
and pick $h_0$ such that 
\begin{equation}\label{paper3:right_ineq}
\tilde{h}<h_0<h_1-N(J_1+J_2).
\end{equation} 
The right inequality of \eqref{paper3:right_ineq} is equivalent to
\[
2(h_0+NJ_2)<2(h_1-NJ_1)
\]
and so we can pick $0<p_1<p_2<1$ such that
\begin{equation*}
2(h_0+NJ_2)<\log(\frac{p_1}{1-p_1})<\log(\frac{p_2}{1-p_2})<2(h_1-NJ_1).
\end{equation*}
By using the last inequalities together with Proposition 4.16 in \cite{Georgii_Haggstrom_Maes_geom}, we can conclude that
\begin{align*}
\mu_{h_0}^{J_2,+}&\leq \gamma_{p_1} \\
\mu_{h_1}^{J_1,+}&\geq \gamma_{p_2}.
\end{align*}
Since $p_1<p_2$ this tells us that $\mu_{h_0}^{J_2,+}\ngeq \mu_{h_1}^{J_1,+}$. On the other hand we have $h_0>\tilde{h}$ which by definition of $\tilde{h}$ implies that $\mu_{h_0}^{J_2,+}\geq \mu_{h_1}^{J_1,+}$. Hence, we get a contradiction and the proof is complete.
\end{proofofprop_no_text}

\subsection{Proof of Theorem \ref{paper3:prop_dom_tree}}\label{paper3:proof_dom_tree}

We will make use of the following lemma from \cite{Liggett_Steif} concerning stochastic domination for completely homogeneous Markov chains on $\T^d$.
\begin{lemma}[Liggett, Steif]\label{paper3:LS3}
Given two 2-state transition matrices $P$ and $Q$, let $\mu_P$ and $\mu_Q$ denote the corresponding completely homogeneous Markov chains on $\T^d$. Then $\mu_P$ dominates $\mu_Q$ if and only if $P(-1,1)\geq Q(-1,1)$ and $P(1,1)\geq Q(1,1)$.
\end{lemma}

\noindent\begin{proofoft}{\ref{paper3:prop_dom_tree}}
To prove \eqref{paper3:prop_dom_tree_ekv1}, let $J_1$, $J_2>0$ and $h_1\in\RR$ be given and note that we get from Lemma~\ref{paper3:LS3} and equation \eqref{paper3:Pt} that $\mu_h^{J_2,+}\geq\mu_{h_1}^{J_1,\pm}$ if and only if 
\[
\frac{e^{t_+(J_2,h)-J_2}}{2\cosh(t_+(J_2,h)-J_2)}\geq \frac{e^{t_\pm(J_1,h_1)-J_1}}{2\cosh(t_\pm(J_1,h_1)-J_1)}
\]
and
\[
\frac{e^{t_+(J_2,h)+J_2}}{2\cosh(t_+(J_2,h)+J_2)}\geq \frac{e^{t_\pm(J_1,h_1)+J_1}}{2\cosh(t_\pm(J_1,h_1)+J_1)}.
\]
Since the map $\RR\ni x \mapsto \frac{e^x}{2\cosh(x)}$ is strictly increasing this is equivalent to
\[
t_+(J_2,h) \geq t_\pm(J_1,h_1) + J_2-J_1 
\]
and
\[
t_+(J_2,h) \geq t_\pm(J_1,h_1) + J_1-J_2 
\]
which in turn is equivalent to
\begin{equation}\label{paper3:cond_stoch_dom_tree}
t_+(J_2,h) \geq t_\pm(J_1,h_1) +|J_1-J_2|=\tau_\pm(J_1,J_2,h_1),
\end{equation}
and so we want to compute the smallest $h\in\RR$ such that \eqref{paper3:cond_stoch_dom_tree} holds. Note that since the map $h\mapsto t_+(J_2,h)$ is strictly increasing and $t_+(J_2,h)\to\pm\infty$ as $h\to\pm\infty$ there always exists such an $h\in\RR$. If $\tau_\pm\geq t^*(J_2)$ or $\tau_\pm<t_-(J_2,-h^*(J_2))$, then the equation
\[
h+d\phi_{J_2}(\tau_\pm)=\tau_\pm
\] 
is equivalent to
\[
t_+(J_2,h)=\tau_\pm
\]
and so in that case the smallest $h\in\RR$ such that \eqref{paper3:cond_stoch_dom_tree} holds is equal to
\[
\tau_\pm-d\phi_{J_2}(\tau_\pm).
\]
If $t_-(J_2,-h^*(J_2))\leq \tau_\pm<t^*(J_2)$, then since $t_+(J_2,h)\geq t^*(J_2)$ whenever $h\geq -h^*(J_2)$ and $t_+(J_2,h)<t_-(J_2,-h^*(J_2))$ whenever $h< -h^*(J_2)$, we have in this case that 
\[
\{\,h\in\RR:\mu_h^{J_2,+}\geq \mu_{h_1}^{J_1,\pm}\,\}=[-h^*(J_2),\infty),
\]
and so the $h$ we are looking for is given by $-h^*(J_2)$.

For \eqref{paper3:prop_dom_tree_ekv2}, we note as above that $\mu_h^{J_2,-}\geq\mu_{h_1}^{J_1,\pm}$ if and only if
\begin{equation}\label{paper3:cond_stoch_dom_tree2}
t_-(J_2,h) \geq \tau_\pm(J_1,J_2,h_1).
\end{equation} 
If $\tau_\pm\leq -t^*(J_2)$ or $\tau_\pm>t_+(J_2,h^*(J_2))$ then we can proceed exactly as in the first case above. If $-t^*(J_2)<\tau_\pm\leq t_+(J_2,h^*(J_2))$, then $t_-(J_2,h)<\tau_\pm$ whenever $h\leq h^*(J_2)$ and $t_-(J_2,h)>\tau_\pm$ whenever $h> h^*(J_2)$ and so in that case we have
\[
\{\,h\in\RR:\mu_h^{J_2,-}\geq \mu_{h_1}^{J_1,\pm}\,\}=(h^*(J_2),\infty),
\]
which yields \eqref{paper3:prop_dom_tree_ekv2} and the proof is complete.
\end{proofoft}

We will now indicate how to compute the bounds in Proposition \ref{paper3:gen_graph} in the special case when $G=\T^d$. By looking at the formula for $f_+$ and using the definition of $h^*$ we get that
\[
f_+(J_1,J_2,h_1)\leq \tau_+(J_1,J_2,h_1)-d\phi_{J_2}(\tau_+(J_1,J_2,h_1)).
\]
Substituting $\tau_+$ and using the bounds $-J\leq \phi_J(t)\leq J$ for all $t\in\RR$ we get the upper bound in Proposition \ref{paper3:gen_graph} with $N=d+1$. For the lower bound, first note that
\[
\begin{split}
\tau_+-d\phi_{J_2}(\tau_+)=h_1+d\big(\phi_{J_1}(t_+(J_1,h_1))&-\phi_{J_2}(t_+(J_1,h_1))\big)+|J_1-J_2| \\
\geq h_1-(d+1)(J_1+J_2).
\end{split}
\]
Moreover it is easy to check that
\[
-h^*(J_2)\geq h_1-(d+1)(J_1+J_2)
\]
when
\[
t_-(J_2,-h^*(J_2))\leq \tau_+\leq t^*(J_2)=t_+(J_2,-h^*(J_2))
\]
and so the lower bound follows at once.
\subsection{Proof of Proposition \ref{paper3:prop_dom_tree2}}

Before we prove anything we would like to recall the fact that we can write (see Remark (ii) after Theorem \ref{paper3:prop_dom_tree})
\[
f_+(J_1,J_2,h_1)=\psi(J_2,\tau_+(J_1,J_2,h_1)) \quad J_1,J_2>0, h_1\in\RR,
\]
where 
\[
\tau_+(J_1,J_2,h_1)=t_+(J_1,h_1)+|J_1-J_2|
\]
and the map $t\mapsto\psi(J_2,t)$ is continuous (see Figure \ref{paper3:fplusfig} for a picture). In the rest of the proof, we will use this fact without further notification. For example, the above immediately gives that $h_1\mapsto f_+(J_1,J_2,h_1)$ is continuous at a point $h_1\in\RR$ if $h_1\mapsto t_+(J_1,h_1)$ is so.

\medskip

\noindent\begin{proofofprop}{\ref{paper3:prop_dom_tree2}}
We will only prove part $a)$ and $c)$. The proof of part $b)$ follows the same type of argument as the proof of part $a)$.

To prove part $a)$, we start to argue that for given $J_1>0$ the map $h_1\mapsto t_+(J_1,h_1)$ is right-continuous at every point $h_1\in\RR$. To see that, take a sequence of reals $\{h_n\}$ such that $h_n\downarrow h_1$ as $n\to\infty$ and note that since the map $h_1\mapsto t_+(J_1,h_1)$ is increasing, the sequence $\{t_+(J_1,h_n)\}$ converges to a limit $\tilde{t}$ with $\tilde{t}\geq t_+(J_1,h_1)$. Moreover, by taking the limit in the fixed point equation we see that
\begin{equation}\label{paper3:tilde}
\tilde{t}=h_1+d\phi_{J_1}(\tilde{t})
\end{equation}
and since $t_+(J_1,h_1)$ is the largest number satisfying \eqref{paper3:tilde} we get $\tilde{t}=t_+(J_1,h_1)$. 

Next, assume $h_1\neq -h^*(J_1)$ and $h_n\uparrow h_1$ as $n\to\infty$. As before, the limit of $\{t_+(J_1,h_n)\}$ exists, denote it by $T$. The number $T$ will again satisfy \eqref{paper3:tilde}. By considering different cases described in Figure \ref{paper3:cases_cont}, we easily conclude that $T=t_+(J_1,h_1)$. Hence, the function $h_1\mapsto t_+(J_1,h_1)$ is continuous for all $h_1\neq -h^*(J)$ and so we get that $h_1\mapsto f_+(J_1,J_2,h_1)$ is also continuous for all $h_1\neq -h^*(J_1)$.
\begin{figure}[!h]
\begin{center}
\scalebox{\scalefactor}{
%
%
\begin{psfrags}%
\psfragscanon%
%
\psfrag{s10}[t][t]{\color[rgb]{0,0,0}\setlength{\tabcolsep}{0pt}\begin{tabular}{c}$t$\end{tabular}}%
\psfrag{s11}[l][l]{\color[rgb]{0,0,0}\setlength{\tabcolsep}{0pt}\begin{tabular}{l}$|h_1|<h^*(J_1)$, $h^*(J_1)>0$\end{tabular}}%
\psfrag{s12}[t][t]{\color[rgb]{0,0,0}\setlength{\tabcolsep}{0pt}\begin{tabular}{c}$t$\end{tabular}}%
\psfrag{s13}[l][l]{\color[rgb]{0,0,0}\setlength{\tabcolsep}{0pt}\begin{tabular}{l}$h_1=h^*(J_1)>0$\end{tabular}}%
\psfrag{s14}[t][t]{\color[rgb]{0,0,0}\setlength{\tabcolsep}{0pt}\begin{tabular}{c}$t$\end{tabular}}%
\psfrag{s15}[l][l]{\color[rgb]{0,0,0}\setlength{\tabcolsep}{0pt}\begin{tabular}{l}$|h_1|>h^*(J_1)$\end{tabular}}%
%
\psfrag{x01}[t][t]{-30}%
\psfrag{x02}[t][t]{-20}%
\psfrag{x03}[t][t]{-10}%
\psfrag{x04}[t][t]{0}%
\psfrag{x05}[t][t]{10}%
\psfrag{x06}[t][t]{20}%
\psfrag{x07}[t][t]{30}%
\psfrag{x08}[t][t]{-30}%
\psfrag{x09}[t][t]{-20}%
\psfrag{x10}[t][t]{-10}%
\psfrag{x11}[t][t]{0}%
\psfrag{x12}[t][t]{10}%
\psfrag{x13}[t][t]{20}%
\psfrag{x14}[t][t]{30}%
\psfrag{x15}[t][t]{-30}%
\psfrag{x16}[t][t]{-20}%
\psfrag{x17}[t][t]{-10}%
\psfrag{x18}[t][t]{0}%
\psfrag{x19}[t][t]{10}%
\psfrag{x20}[t][t]{20}%
\psfrag{x21}[t][t]{30}%
%
\psfrag{v01}[r][r]{-20}%
\psfrag{v02}[r][r]{0}%
\psfrag{v03}[r][r]{20}%
\psfrag{v04}[r][r]{-20}%
\psfrag{v05}[r][r]{0}%
\psfrag{v06}[r][r]{20}%
\psfrag{v07}[r][r]{-20}%
\psfrag{v08}[r][r]{0}%
\psfrag{v09}[r][r]{20}%
%
\resizebox{12cm}{!}{\includegraphics{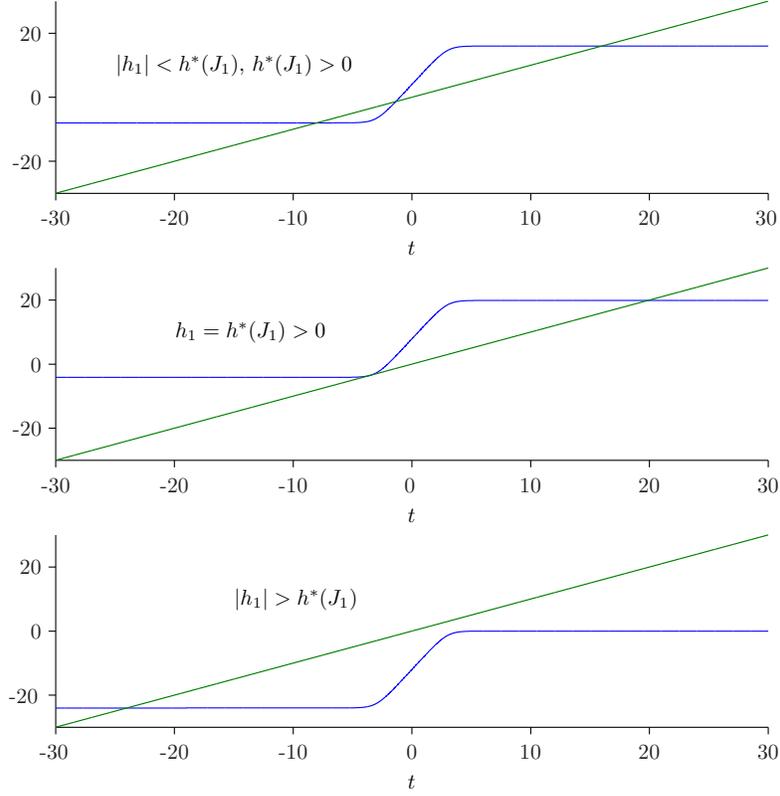}}%
\end{psfrags}%
%
}
\caption{A picture of the different cases in the fixed point equation that can occur when $h_1\neq -h^*(J_1)$. Here, $d=4$ and $J_1=3$.}
\label{paper3:cases_cont}
\end{center}
\end{figure}

Now assume $h_1=-h^*(J_1)$. By considering sequences $h_n\downarrow -h^*(J_1)$ and $h_n\uparrow -h^*(J_1)$ we can similarly as above see that
\begin{align*}
\tau_+(J_1,J_2,-h^*(J_1)+):&=\lim_{h\downarrow -h^*(J_1)}\tau_+(J_1,J_2,h)=t_+(J_1,-h^*(J_1))+|J_1-J_2| \\
\tau_+(J_1,J_2,-h^*(J_1)-):&=\lim_{h\uparrow -h^*(J_1)}\tau_+(J_1,J_2,h)=t_-(J_1,-h^*(J_1))+|J_1-J_2|
\end{align*}
and so 
\[
\tau_+(J_1,J_2,-h^*(J_1)+)=\tau_+(J_1,J_2,-h^*(J_1)-) \quad \iff \quad h^*(J_1)=0.
\]
Since $h^*(J_1)=0$ if and only if $0<J_1\leq J_c$ the continuity of $h_1\mapsto f_+(J_1,J_2,h_1)$ at $-h^*(J_1)$ follows at once in that case. If $J_1=J_2$, then
\begin{align*}
\tau_+(J_1,J_2,-h^*(J_1)+)&=t_+(J_2,-h^*(J_2)) \\
\tau_+(J_1,J_2,-h^*(J_1)-)&=t_-(J_2,-h^*(J_2))
\end{align*}
and since 
\[
\psi(J_2,t_+(J_2,-h^*(J_2)))=\psi(J_2,t_-(J_2,-h^*(J_2))),
\]
the continuity is clear also in that case. If $J_1>J_c$ and $0<J_2\leq J_c$, then 
\[
\tau_+(J_1,J_2,-h^*(J_1)+)\neq\tau_+(J_1,J_2,-h^*(J_1)-) 
\]
and the map $t\mapsto \psi(J_2,t)$ becomes strictly increasing, hence $h_1\mapsto f_+(J_1,J_2,h_1)$ is discontinuous at $-h^*(J_1)$. For the case when $J_1>J_c$, $J_2>J_c$, $J_1\neq J_2$ just note that $h_1\mapsto f_+(J_1,J_2,h_1)$ is continuous at $-h^*(J_1)$ if and only if $a$ and $b$ (defined in the statement of the proposition) are in the flat region in the upper graph of Figure \ref{paper3:fplusfig}.

To prove part $c)$ we take a closer look at the map $(J_2,t)\mapsto\psi(J_2,t)$. By definition, this map is
\[
\psi(J_2,t)=\begin{cases}
-h^*(J_2) &\text{if} \quad t_-(J_2,-h^*(J_2))\leq t<t^*(J_2) \\
t-d\dpt\phi_{J_2}(t) &\text{if} \quad t\geq t^*(J_2) \text{ or } t<t_-(J_2,-h^*(J_2)).
\end{cases}
\]
From the continuity of $t\mapsto \psi(J_2,t)$ for fixed $J_2$ and the facts that $J_2\mapsto t^*(J_2)$, $J_2\mapsto t_-(J_2,-h^*(J_2))$, $J_2\mapsto -h^*(J_2)$ and $(J_2,t)\mapsto t-d\phi_{J_2}(t)$ are all continuous, we get that $\psi$ is (jointly) continuous and so the result follows. 
\end{proofofprop}

\subsection{Proof of Proposition \ref{paper3:prop_dom_tree3}}

\begin{proofofprop_no_text}
To prove the statement, we will show that the inequality
\begin{equation}\label{paper3:der_tplus}
\frac{\partial}{\partial J}t_+(J,h)\geq 1
\end{equation}
holds if $a)$ $h\geq 0$ and $J\geq J_c$ or $b)$ $h<0$ and $h^*(J)>-h$.
By integrating equation \eqref{paper3:der_tplus} the statement follows. The proof of equation \eqref{paper3:der_tplus} will be an easy modification of the proof of Lemma 5.2 in \cite{Liggett_Steif}. The proof is quite short and so we give a full proof here, even though it is more or less the same as the proof in \cite{Liggett_Steif}.

Write $\phi(J,t)$ for $\phi_J(t)$ and use subscripts to denote partial derivatives. By differentiating the relation
\[
h+d\,\phi(J,t_+(J,h))=t_+(J,h)
\]
with respect to $J$ and solving, we get
\[
\frac{\partial}{\partial J}t_+(J,h)=\frac{d\,\phi_1(J,t_+(J,h))}{1-d\,\phi_2(J,t_+(J,h))}.
\]
To get the left hand side bigger or equal to one, we need
\begin{equation}\label{paper3:ineq1}
d\,\phi_2(J,t_+(J,h))<1
\end{equation} 
and
\begin{equation}\label{paper3:ineq2}
\phi_1(J,t_+(J,h))+\phi_2(J,t_+(J,h))\geq \frac{1}{d}.
\end{equation}
The first inequality is immediate since in the cases $a)$ and $b)$ above, the function $t\mapsto h+d\,\phi(J,t)$ crosses the line $t\mapsto t$ from above to below. For \eqref{paper3:ineq2}, note that 
\begin{align*}
\phi_1(J,t)&=\frac{1}{2}\big(\tanh(J+t)-\tanh(J-t)\big) \\
\phi_2(J,t)&=\frac{1}{2}\big(\tanh(J+t)+\tanh(J-t)\big)
\end{align*}
and so
\[
\phi_1(J,t)+\phi_2(J,t)=\tanh(J+t),
\]
which yields that $\phi_1+\phi_2$ is increasing in both variables. Moreover, since $\tanh(J_c)=\frac{1}{d}$ (see \cite{Georgii}), we get
\[
\phi_1(J_c,0)+\phi_2(J_c,0)=\frac{1}{d}
\]
and so 
\begin{equation}\label{paper3:ineq3}
\phi_1(J,t)+\phi_2(J,t)\geq \frac{1}{d} \quad \text{if $J\geq J_c$, $t\geq 0$.}
\end{equation}
To complete the proof, observe that in the cases $a)$ and $b)$, we have $J\geq J_c$ and $t_+(J,h)\geq 0$.
\end{proofofprop_no_text}

\subsection{Proof of Proposition \ref{paper3:fuzzy_prod_Zd}}

In the proof we will use the following results from \cite{Liggett_Steif} concerning domination of product measures.
\begin{defin}[Downward FKG, Liggett, Steif]
Given a finite or countable set $V$, a measure $\mu$ on $\{-1,1\}^V$ is called downward FKG if for any finite $A\subseteq V$, the conditional measure $\mu(\,\cdot\, |\,\eta\equiv 0\text{ on }A\,)$ has positive correlations.
\end{defin}
Here, as usual, positive correlations is defined as follows:
\begin{defin}[Positive correlations]
A probability measure $\mu$ on $\{-1,1\}^{V}$ where $V$ is a finite or countable set is said to have positive correlations if
\[
\int fg\,d\mu \geq \int f\,d\mu\int g\,d\mu
\]
for all real-valued, continuous and increasing functions $f,g$ on $\{-1,1\}^{V}$.
\end{defin}
\begin{theorem}[Liggett, Steif]\label{paper3:LS4}
Let $\mu$ be a translation invariant measure on \newline $\{-1,1\}^{\Zd}$ which also is downward FKG and let $p\in [0,1]$. Then the following are equivalent:
\begin{itemize}
\item[a)] $\mu\geq\gamma_p$. 
\item[b)] $\displaystyle\limsup_{n\to\infty}\dpt\mu(\,\eta\equiv -1 \textrm{ on } [1,n]^d\,)^{1/n^d}\leq 1-p$.   
\end{itemize}
\end{theorem}

\noindent
\Remarks 
\begin{itemize}
\item[(i)] In particular, Theorem \ref{paper3:LS4} gives us that if two translation invariant, downward FKG measures have the same above limsup, then they dominate the same set of product measures.
\item[(ii)] In \cite{Liggett_Steif} they had a third condition in Theorem \ref{paper3:LS4} which we will not use and so we simply omit it.
\end{itemize}
Before we state the next lemma we need to recall the following definition. 
\begin{defin}[FKG lattice condition]
Suppose $V$ is a finite set and let $\mu$ be a probability measure on $\{-1,1\}^V$ which assigns positive probabilty to each element. For $\eta$, $\xi\in\{-1,1\}^V$ define $\eta\vee\xi$ and $\eta\wedge\xi$ by
\[
(\eta\vee\xi)(x)=\max(\eta(x),\xi(x)),\,(\eta\wedge\xi)(x)=\min(\eta(x),\xi(x)), \,x\in V. 
\]
We say that $\mu$ satisfies the FKG lattice condition if 
\[
\mu(\eta\wedge\xi)\mu(\eta\vee\xi)\geq\mu(\eta)\mu(\xi)
\]
for all $\eta$, $\xi\in\{-1,1\}^V$
\end{defin}
Given a measure $\mu$ on $\{-1,1\}^{\Zd}$ we will denote its projection on $\{-1,1\}^T$ for finite $T\subseteq\Zd$ by $\mu_T$.
\begin{lemma}\label{paper3:FKG_fuzzy}
The measures $\nu_{q,J,r}^{\Zd,\pm}$ are FKG in the sense that $\nu_{T,q,J,r}^{\Zd,\pm}$ satisfies the FKG lattice condtion for each finite  $T\subseteq\Zd$.
\end{lemma}
\begin{proof}
For $n\geq 2$, let $\Lambda_n=\{-n,\dots,n\}^d$ and denote the finite volume Potts measures on $\{-1,1\}^{\Lambda_n}$ with boundary condition $\eta\equiv 1$ and  $\eta\equiv q$ by $\pi_{q,J}^{n,1}$ and $\pi_{q,J}^{n,q}$. Furthermore, let $\nu_{q,J,r}^{n,-}$ and $\nu_{q,J,r}^{n,+}$ denote the corresponding fuzzy Potts measures. Given the convergence in the Potts model, it is clear that $\nu_{T,q,J,r}^{n,\pm}$ converges weakly to $\nu_{T,q,J,r}^{\Zd,\pm}$ as $n\to\infty$ for each finite $T\subseteq \Zd$. Since the FKG lattice condition is closed under taking projections (see \cite[p.~28]{Grimmett_RC_model}) and weak limits we are done if we can show that $\nu_{q,J,r}^{n,\pm}$ satisfies the FKG lattice condition for each $n\geq 2$. In \cite{Haggstrom_pos_corr} it is proved that for an arbitrary finite graph $G=(V,E)$ the finite volume fuzzy Potts measure with free boundary condition and parameters $q$, $J$, $r$ is monotone in the sense that 
\begin{equation}\label{paper3:proof:monotone}
\nu_{q,J,r}^G(Y(x)=1\,|\,Y(V\setminus\{x\})=\eta)\leq\nu_{q,J,r}^G(Y(x)=1\,|\,Y(V\setminus\{x\})=\eta^\prime)
\end{equation}
for all $x\in V$ and $\eta$, $\eta^\prime\in\{-1,1\}^{V\setminus\{x\}}$ with $\eta\leq\eta^\prime$. We claim that it is possible to modify the argument given there to prove that $\nu_{q,J,r}^{n,\pm}$ are monotone for each $n\geq 2$. (Recall from \cite{Grimmett_RC_model} the fact that if $V$ is finite and $\mu$ is a probabilty measure on $\{-1,1\}^V$ that assigns positive probabilty to each element, then monotone is equivalent to the FKG lattice condition.) The proof of \eqref{paper3:proof:monotone} is quite involved. However, the changes needed to prove our claim are quite straightforward and so we will only give an outline for how that can be done. Furthermore, we will only consider the minus case, the plus case is similar.
 
By considering a sequence $\eta=\eta_1\leq\eta_2\leq\dots\leq\eta_m=\eta^\prime$ where $\eta_i$ and $\eta_{i+1}$ differ only at a single vertex, it is easy to see that it is enough to prove that for all $x$, $y\in\Lambda_n$ and $\eta\in\{-1,1\}^{\Lambda_n\setminus\{x,y\}}$ we have
\begin{equation}\label{paper3:proof:cond_indep}
\begin{split}
&\nu_{q,J,r}^{n,-}(Y(x)=1,Y(y)=1\,|\,Y(\Lambda_n\setminus\{x,y\})=\eta) \\
&\quad\geq\nu_{q,J,r}^{n,-}(Y(x)=1\,|\,Y(\Lambda_n\setminus\{x,y\})=\eta) \\
&\qquad\cdot\nu_{q,J,r}^{n,-}(Y(y)=1\,|\,Y(\Lambda_n\setminus\{x,y\})=\eta).
\end{split}
\end{equation} 
Fix $n\geq 2$, $x$, $y$ and $\eta$ as above. We will first consider the case when $x$ and $y$ are not neighbors. At the end we will see how to modify the argument to work when $x$, $y$ are neighbors as well. Define $V_-=\{z\in\Lambda_n\setminus\{x,y\}:\eta(z)=-1\}$ and $V_+=\{z\in\Lambda_n\setminus\{x,y\}:\eta(z)=1\}$. Furthermore, denote by $E_n$ the set of edges $\langle u,v\rangle$ with either $u$, $v\in\Lambda_n$ or $u\in\Lambda_n$, $v\in\partial\Lambda_n$ and let $\P$ denote the probability measure on $W=\{1,\dots,q\}^{\Lambda_n\cup\partial \Lambda_n}\times\{0,1\}^{E_n}$ which to each site $u\in\Lambda_n\cup\partial\Lambda_n$ chooses a spin value uniformly from $\{1,\dots,q\}$, to each edge $\langle u,v\rangle$ assigns value $1$ or $0$ with probabilities $p$ and $1-p$ respectively and which does those things independently for all sites and edges. Define the following events on $W$
\begin{align*}
A&=\{(\sigma,\xi):\,(\sigma(u)-\sigma(v))\xi(e)=0,\,\forall e=\langle u,v\rangle\in E_n\,\}, \\
B&=\{(\sigma,\xi):\,\sigma(z)\in\{1,\dots,r\}\,\forall z\in V_-,\,\sigma(z)\in\{r+1,\dots,q\}\,\forall z\in V_+\}, \\
C&=\{(\sigma,\xi):\sigma(z)=1,\,\forall z\in\partial\Lambda_n\,\},
\end{align*}
and let $\P^\prime$ and $\P^{\prime\prime}$ be the probability measures on $\{1,\dots,q\}^{\Lambda_n}\times\{0,1\}^{E_n}$ obtained from $\P$ by conditioning on $A\cap C$ and $A\cap B\cap C$ respectively. ($\P^\prime$ is usually referred to as the Edward-Sokal coupling, see \cite{Georgii_Haggstrom_Maes_geom}.) It is well known (and easy to check) that the spin marginal of $\P^\prime$ is $\pi_{q,J}^{n,1}$ and that the edge marginal is the so called random-cluster measure defined as the probability measure on $\{0,1\}^{E_n}$ which to each $\xi\in\{0,1\}^{E_n}$ assigns probability proportional to 
\[
q^{k(\xi)}\prod_{e\in E_n}p^{\xi(e)}(1-p)^{1-\xi(e)},
\]
where $k(\xi)$ is the number of connected components in $\xi$ not reaching $\partial \Lambda_n$. In a similar way it is possible (by counting) to compute the spin and edge marginal of $\P^{\prime\prime}$: The spin marginal $\pi^{\prime\prime}$ is simply $\pi_{q,J}^{n,1}$ conditioned on $B$ and the edge marginal $\phi^{\prime\prime}$ assigns probability to a configuration $\xi\in\{0,1\}^{E_n}$ proportional to  
\[
1_D r^{k_0(\xi)}(q-r)^{k_1(\xi)}q^{k_x(\xi)+k_y(\xi)}\prod_{e\in E_n}p^{\xi(e)}(1-p)^{1-\xi(e)},
\]
where $k_0(\xi)$ is the number of clusters intersecting $V_-$ but not reaching $\partial \Lambda_n$, $k_1(\xi)$ is the number of clusters intersecting $V_+$, $k_x(\xi)$ (resp $k_y(\xi)$) is $1$ if $x$ (resp $y$) is in a singleton connected component and $0$ otherwise and $D$ is the event that no connected component in $\xi$ intersects both $V_-$ and $V_+$. Observe that \eqref{paper3:proof:cond_indep} is the same as
\begin{equation}\label{paper3:proof:goal}
\begin{split}
&\pi^{\prime\prime}(X(x)\in\{r+1,\dots,q\},X(y)\in\{r+1,\dots,q\}) \\
&\quad\geq\pi^{\prime\prime}(X(x)\in\{r+1,\dots,q\})\,\pi^{\prime\prime}(X(y)\in\{r+1,\dots,q\}).
\end{split}
\end{equation}
An important feature of the coupling $\P^{\prime\prime}$ is that it gives a way to obtain a spin configuration $X\in\{1,\dots,q\}^{\Lambda_n}$ distributed as $\pi^{\prime\prime}$: 
\begin{enumerate}
\item Pick an edge configuration $\xi$ according to $\phi^{\prime\prime}$. 
\item Assign $X=1$ to the connected components of $\xi$ that intersect $\partial \Lambda_n$ and denote the union of those components by $\tilde{C}$.
\item Assign independently spins to a connected component $C\neq\tilde{C}$ of $\xi$ where the spin is taken according to the uniform distribution on
\begin{center}
\begin{tabular}{ll}
$\{1,\dots,r\}$ & if $C$ intersects $V_-$, \\
$\{r+1,\dots,q\}$ & if $C$ intersects $V_+$, \\
$\{1,\dots,q\}$ & if $C$ is a singleton vertex $x$ or $y$.  
\end{tabular}
\end{center}
\end{enumerate}
By defining the functions $f_x$, $f_y:\{0,1\}^{E_n}\to\RR$ as
\[
f_x(\xi)=\begin{cases}
0, \text{ if $C_x=\tilde{C}$ or $C_x$ intersects $V_-$,} \\
\frac{q-r}{q}, \text{ if $C_x$ is a singleton,} \\
1, \text{ otherwise,}
\end{cases}
\]
where $C_x$ is the connected component of $\xi$ containing $x$ ($f_y$ defined analogously), we see as in \cite{Haggstrom_pos_corr} that \eqref{paper3:proof:goal} follows if 
\begin{equation}\label{paper3:proof:goal2}
\int f_xf_y\,d\phi^{\prime\prime} \geq \int f_x\,d\phi^{\prime\prime} \int f_y\,d\phi^{\prime\prime}.
\end{equation}
The significance of $f_x$ and $f_y$ is that $f_x(\xi)$ is the conditional probability that $X(x)\in\{r+1,\dots,q\}$ given $\xi$ and similarly for $f_y$, and that the events  $X(x)\in\{r+1,\dots,q\}$ and $X(y)\in\{r+1,\dots,q\}$ are conditionally independent given $\xi$. With all this setup done it is a simple task to see that to prove \eqref{paper3:proof:goal2} we can proceed exactly as in \cite[p.~1154-1155]{Haggstrom_pos_corr}. 

To take care of the case when $x$ and $y$ are neighbors, observe that everything we have done so far also works for the graph with one edge deleted, i.e.\ the graph with vertex set $\Lambda_n$ and edge set $E_n\setminus\{\langle x,y\rangle\}$. Hence we can get \eqref{paper3:proof:goal} for that graph. However the observation in \cite[1156]{Haggstrom_pos_corr} gives us \eqref{paper3:proof:goal} even in the case when we reinsert the edge $\langle x,y\rangle$. 
\end{proof}

\noindent
\begin{proofofprop}{\ref{paper3:fuzzy_prod_Zd}}
Let $k$, $l\in\{0,-,+\}$ be given and let $A_n=[1,n]^d$, $n\geq 2$. We are done if there exists $0<c<1$ (independent of $k$, $l$ and $n$) such that
\[
\nu_{q,J,r}^{\Zd,k}(\,\eta\equiv -1\text{ on }A_n\,)\geq c^{|\partial A_n|}\nu_{q,J,r}^{\Zd,l}(\,\eta\equiv -1\text{ on }A_n\,) \text{ for all $n$.}
\]
As for the Ising model, it is well known  that the infinite volume Potts measures satisfy the so called uniform nonnull property (sometimes called uniform finite energy property), which means that for some $c>0$, the conditional probability of having a certain spin at a given site given everything else is at least $c$. (See for example \cite{Haggstrom_Gibsian_Zd} for a more precise definition.) We get for arbitrary $\sigma\in\{1,\dots,q\}^{\partial A_n}$ 
\begin{equation}\label{paper3:ineq_fuzzy_Zd}
\nu_{q,J,r}^{\Zd,k}(\,\eta\equiv -1\text{ on }[1,n]^d\,)\geq c^{|\partial A_n|} \pi_{q,J}^{A_n,\sigma}(\,Y\equiv -1  \text{ on $A_n$}\,).
\end{equation}
Since $\nu_{q,J,r}^{\Zd,l}(\,\eta\equiv -1\text{ on }[1,n]^d\,)$ can be written as a convex combination of the terms in the far right side of \eqref{paper3:ineq_fuzzy_Zd} the result follows at once.
\end{proofofprop}

\subsection{Proof of Proposition \ref{paper3:fuzzy_prod_Td}}
Let $\rho$ denote the root of $\Td$ and let $V_n$ be the set of all sites in $\Td$ with distance at most $n$ from $\rho$. If $x$ is on the unique self-avoiding path from $\rho$ to $y$, we say that $y$ is a descendant of $x$. Given $x\in\Td$, let $S_x$ denote the set of vertices of all descendants of $x$ (including $x$). Moreover, let $T_x$ denote the subtree of $\Td$ whose vertex set is $S_x$ and edge set consisting of all edges $\langle u,v \rangle \in E(\Td)$ with $u$, $v\in S_x$. In the proof of Proposition \ref{paper3:fuzzy_prod_Td}, we will use the following Lemma from \cite{Liggett_Steif}:   
\begin{prop}[Liggett, Steif]\label{paper3:LS5}
Let $p\in [0,1]$, $\{\,P(i,j):\,i,j\in\{\ -1,1\}\,\}$ be a transition matrix for an irreducible 2-state Markov chain with $P(-1,1)\leq P(1,1)$ and let $\mu$ be the distribution of the corresponding completely homogeneous Markov chain on $\Td$. Then the following are equivalent:
\begin{itemize}
\item[a)] $\mu \geq\gamma_p$.
\item[b)] $\displaystyle\limsup_{n\to\infty}\mu(\,\eta\equiv -1 \text{ on $V_n$}\,)^{1/|V_n|}\leq 1-p$.
\item[c)] $P(-1,1)\geq p$.
\end{itemize}
\end{prop}

\noindent\begin{proofofprop}{\ref{paper3:fuzzy_prod_Td}}
Fix $J>0$, $q\geq 3$ and $r\in\{1,\dots,q-1\}$ with $e^{2J}\geq q-2$. In \cite{Haggstrom_Gibsian_tree}, it is proved that $\nu_{q,J,r}^{\Td,0}$ is a completely homogeneous Markov chain on $\Td$ for all values of the parameters with transition matrix
\begin{equation*}\label{paper3:trans_fuzzy}
\begin{pmatrix}
\frac{e^{2J}+r-1}{e^{2J}+q-1}& \frac{q-r}{e^{2J}+q-1}\\
\frac{r}{e^{2J}+q-1}& \frac{e^{2J}+q-r-1}{e^{2J}+q-1}\\
\end{pmatrix}.
\end{equation*}
Hence, from Proposition \ref{paper3:LS5} we get that $\nu_{q,J,r}^{\Td,0}\geq \gamma_p$ if and only if
\begin{equation}\label{paper3:dom_free}
p\leq \frac{q-r}{e^{2J}+q-1}.
\end{equation}
Furthermore, in \cite[p.~10]{Haggstrom_Gibsian_tree} the authors also 
derive the transition matrix for $\pi_{q,J}^{\Td,1}$ from which we can compute the following:
\[
\begin{split}
\nu_{q,J,r}^{\Td,-}(\,\eta\equiv -1 \text{ on $V_n$}\,)&\geq\displaystyle\sum_{i=1}^r\pi_{q,J}^{\Td,1}(\,X\equiv i \text{ on $V_n$}\,) \\
&=\frac{b}{b+q-1}\left(\frac{ce^{2J}}{ce^{2J}+q-1}\right)^{|V_n|-1} \\
&\quad+\frac{r-1}{b+q-1}\left(\frac{e^{2J}}{c+e^{2J}+q-2}\right)^{|V_n|-1}
\end{split}
\]
where
\begin{align*}
b&=\frac{\pi_{q,J}^{\Td,1}(\,X(\rho)=1\,)}{\pi_{q,J}^{\Td,1}(\,X(\rho)=2\,)} \\
c&=\frac{\pi_{q,J}^{T_x,1}(\,X(x)=1\,)}{\pi_{q,J}^{T_x,1}(\,X(x)=2\,)}, \quad x\neq\rho.
\end{align*}
(Of course, homogeneity gives that the last quotient is independent of $x$.) We get that
\begin{equation}\label{paper3:limsup_r1}
\begin{split}
&\limsup_{n\to\infty} \nu_{q,J,r}^{\Td,-}(\,\eta\equiv -1 \text{ on $V_n$}\,)^{1/|V_n|} \\
&\quad\geq \frac{ce^{2J}}{ce^{2J}+q-1}+\frac{e^{2J}}{c+e^{2J}+q-2}.
\end{split}
\end{equation}
Now, assume that the underlying Gibbs measures for the Potts model satisfy $\pi^{\Td,1}_{q,J}\neq \pi^{\Td,0}_{q,J}$. It is known \cite{Aizenman_et_al_conv_Potts} that this is equivalent to having
\[
\pi^{\Td,1}_{q,J}(\,X(x)=1\,)>\frac{1}{q},\quad\forall\, x\in\Td.
\]
In \cite{Haggstrom_Gibsian_tree}, the authors observed that if $a=\pi_{q,J}^{\Td,1}(\,X(\rho)=1\,)$, then from symmetry reasons
\[
b=\frac{(q-1)a}{1-a}.
\]
Hence, if $a>\frac{1}{q}$ we get $b>1$. Moreover, from the recursion formula in \cite[p.~9]{Haggstrom_Gibsian_tree} we obtain
\begin{equation}\label{paper3:rec}
b=\frac{(ce^{2J}+q-1)^{d+1}}{(c+e^{2J}+q-2)^{d+1}}.
\end{equation}
It is easy to see from \eqref{paper3:rec} that if $b>1$ then $c>1$. Hence, we can choose $p\in (0,1)$ such that
\begin{equation}\label{paper3:choose_p}
\frac{q-r}{ce^{2J}+q-1}< p\leq \frac{q-r}{e^{2J}+q-1}.
\end{equation}
Moreover, an easy calculation gives us that
\[
\frac{ce^{2J}}{ce^{2J}+q-1}+\frac{e^{2J}}{c+e^{2J}+q-2}\geq \frac{ce^{2J}+q-2}{ce^{2J}+q-1} 
\]
and since 
\[
1-p<\frac{ce^{2J}+r-1}{ce^{2J}+q-1}\leq\frac{ce^{2J}+q-2}{ce^{2J}+q-1}
\]
we get from \eqref{paper3:limsup_r1}
\[
\limsup_{n\to\infty}\nu_{q,J,r}^{\Td,-}(\,\eta\equiv -1 \text{ on $V_n$}\,)^{1/|V_n|}>1-p.
\]
It is now clear that for $p$ as in \eqref{paper3:choose_p} we have that $\nu_{q,J,r}^{\Td,0}$ dominates $\gamma_p$ but $\nu_{q,J,r}^{\Td,-}$ does not dominate $\gamma_p$.
\end{proofofprop}

\noindent
\Remark By deriving the transition matrix for $\pi_{q,J}^{\Td,q}$ it is probably possible to prove that there exists $p\in(0,1)$ such that $\nu_{q,J,r}^{\Td,0}$ dominates $\gamma_p$ but $\nu_{q,J,r}^{\Td,+}$ does not dominate $\gamma_p$. 

\section{Conjectures}

We end with the following conjectures concerning the fuzzy Potts model. The corresponding statements for the Ising model are proved in \cite{Liggett_Steif}.

\begin{conjecture}\label{paper3:conj_stoch_dom_Zd}
Let $q\geq 3$, $r\in \{1,\dots,q-1\}$ and consider the fuzzy Potts model on $\Zd$. If $J_1$, $J_2>0$ with $J_1\neq J_2$, then $\nu_{q,J_1,r}^{\Zd,+}$ and $\nu_{q,J_2,r}^{\Zd,+}$ are not stochastically ordered.
\end{conjecture}
\begin{conjecture}\label{paper3:conj_stoch_dom_Zd_sup_prod}
Let $q\geq 3$, $r\in \{1,\dots,q-1\}$ and consider the fuzzy Potts model on $\Zd$. If $0<J_1< J_2$, then 
\[
\sup\{\,p\in[0,1]:\,\nu_{q,J_1,r}^{\Zd,+}\geq\gamma_p\,\}>\sup\{\,p\in[0,1]:\,\nu_{q,J_2,r}^{\Zd,+}\geq\gamma_p\,\}.
\]
\end{conjecture}
\begin{conjecture}\label{paper3:conj_stoch_dom_Td_prod}
Let $J>0$, $q\geq 3$, $r\in \{1,\dots,q-1\}$ and consider the fuzzy Potts model on $\Td$. Define the sets:
\begin{equation}\label{paper3:conj_sets}
\begin{split}
D_+&=\{\,p\in[0,1]:\,\nu_{q,J,r}^{\Zd,+}\geq\gamma_p\,\}, \\
D_-&=\{\,p\in[0,1]:\,\nu_{q,J,r}^{\Zd,-}\geq\gamma_p\,\}, \\
D_0&=\{\,p\in[0,1]:\,\nu_{q,J,r}^{\Zd,0}\geq\gamma_p\,\},
\end{split}
\end{equation}
If the underlying Gibbs measures for the Potts model satisfy $\pi_{q,J}^{\Td,1}\neq \pi_{q,J}^{\Td,0}$, then the sets in \eqref{paper3:conj_sets} are all different from each other.
\end{conjecture}
\begin{conjecture}\label{paper3:conj_stoch_dom_Td}
Let $q\geq 3$, $r\in \{1,\dots,q-1\}$ and consider the fuzzy Potts model on $\Td$. Denote the critical value corresponding to non-uniqueness of Gibbs states for the Potts model by $J_c$. If $J_c<J_1<J_2$ then $\nu_{q,J_1,r}^{\Td,+}\leq\nu_{q,J_2,r}^{\Td,+}$.
\end{conjecture}
\noindent\Remark If $J_1<J_2< J_c$, then 
\[
\nu_{q,J_1,r}^{\Td,+}(\,\eta(x)=1\,)=\nu_{q,J_2,r}^{\Td,+}(\,\eta(x)=1\,)=\frac{q-r}{q}
\]
and so in that case, $\nu_{q,J_1,r}^{\Td,+}$ and $\nu_{q,J_2,r}^{\Td,+}$ can not be stochastically ordered.

\section*{Acknowledgement}

The author wants to thank Jeffrey Steif for presenting the problems, a careful reading of the manuscript and for valuable comments and discussions.


\bibliographystyle{amsplain}


\end{document}